\documentclass[12pt]{extarticle}
\usepackage{amsmath, amsthm, amssymb}
\usepackage[usenames,dvipsnames,table]{xcolor}
\usepackage[colorlinks=true,linkcolor=blue,urlcolor=blue]{hyperref}
\usepackage{graphicx}
\usepackage{caption}
\usepackage{mathtools}
\usepackage{enumerate}
\usepackage{multirow}
\usepackage{verbatim}
\usepackage{tikz,tikz-cd,tikz-3dplot}
\usepackage{amssymb}
\usetikzlibrary{matrix}
\usetikzlibrary{arrows} 

\newtheorem{theorem}{Theorem}
\newtheorem*{theorem*}{Theorem}
\newtheorem{proposition}[theorem]{Proposition}
\newtheorem{lemma}[theorem]{Lemma}

\theoremstyle{definition}
\newtheorem{definition}[theorem]{Definition}

\newtheorem{remark}[theorem]{Remark}

\newtheorem{example}[theorem]{Example}

\numberwithin{theorem}{section}
\usepackage[noend]{algpseudocode}
\usepackage{caption}
\usepackage[normalem]{ulem}
\usepackage{subcaption}
\oddsidemargin \evensidemargin
\usepackage{makecell}
\usepackage{array}
\usepackage[giveninits=true, maxbibnames=99, style=alphabetic,url=false]{biblatex}
\usepackage{amsfonts} 
\usepackage{mathtools}
\usepackage{wrapfig}
\usepackage{cleveref}
\usepackage{amssymb}
\usepackage{tikz}
\usepackage{tikz-cd}
\usepackage[affil-it]{authblk}
\usepackage{blkarray}
\usepackage[english]{babel}
\usepackage[ruled,vlined]{algorithm2e}

\def\kk{\mathbb{K}}
\def\oK{{\overline \kk}}

\def\Ac{\mathcal{A}}		

\def\PP{\mathbb{P}}

\def\RR{{\mathbb R}}
\def\ZZ{{\mathbb Z}}

\def\f{\mathbf{f}}

\def\u{\mathbf{u}}
\def\x{\mathbf{x}}

\DeclareMathOperator{\Secpol}{\Sigma-pol}

\DeclareMathOperator{\Hom}{Hom}

\DeclareMathSymbol{\ast}{\mathbin}{symbols}{"03}

\DeclareMathOperator{\Res}{Res}
\DeclareMathOperator{\ResA}{Res_{\mathcal{A}}}
\DeclareMathOperator{\X}{\mathbf{X}}

\DeclareMathOperator{\F}{\mathbf{F}}
\DeclareMathOperator{\G}{\mathbf{G}}
\DeclareMathOperator{\D}{\mathcal{D}}

\DeclareMathOperator{\Elim}{Elim}

\DeclareMathOperator{\Vol}{Vol}
\DeclareMathOperator{\codim}{codim}

\DeclareMathOperator{\Nef}{Nef}

\DeclareMathOperator{\init}{init}

\DeclareMathOperator{\spann}{span}
\DeclareMathOperator{\conv}{conv}

\DeclareMathOperator{\rc}{rc}

\DeclareMathOperator{\A}{A}

\DeclareMathOperator{\Cay}{Cay}

\DeclareMathOperator{\MV}{MV}

\DeclareMathOperator{\rank}{rank}

\DeclareMathOperator{\Cl}{Cl}

\newcommand{\cc}[1]{\textcolor{blue}{#1}}

\title{\bf A Sparse Overview on Sparse Resultants}
\author{Carles Checa, Ioannis Z. Emiris, Christos Konaxis}

\date{}

\begin{document}
\maketitle

\begin{abstract}
  \noindent

In this survey, we give an overview of advances in the theory and computation of sparse resultants. First, we examine the construction and proof of the Canny-Emiris formula, which gives a rational determinantal formula. Second, we discuss and compare the latter with the computation of the sparse resultant as the determinant of the Koszul complex given by $n + 1$ nef divisors in a toric variety. Finally, we cover techniques for computing the Newton polytope of sparse resultants.

\end{abstract}

\section{Introduction}\label{sec:intro}

Resultants provide a key symbolic technique for eliminating variables and solving systems of polynomial equations. In many cases, it is the preferred method because it can eliminate multiple variables simultaneously, handle parametric coefficients, and efficiently solve specific types of zero-dimensional systems by reducing the question to numerical linear algeba.
The \emph{resultant} of an overconstrained polynomial system characterizes the existence of common roots through algebraic conditions on the input coefficients.

Various methods exist for constructing \emph{resultant matrices}, i.e.\ matrices whose determinant equals the resultant or a nontrivial multiple of it.
These matrices are valuable both for computing the resultant polynomial itself and, more importantly, for solving systems of polynomial equations \cite[Chapter~3]{coxlitosh}.
Specific instances include the high-school example of the coefficient matrix of $n+1$ linear polynomials, the Sylvester matrix of a polynomial pair, and the Bezoutian.

\begin{figure}[t]
    \centering
    \includegraphics[width=0.75\textwidth]{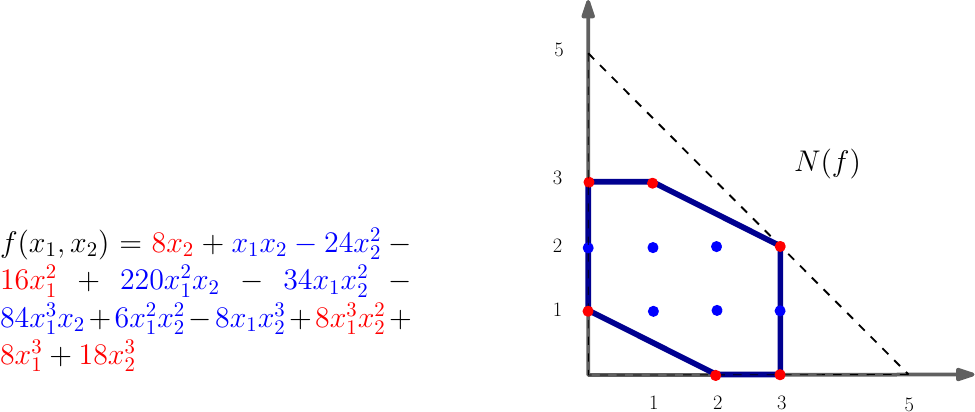}
    \caption{Newton polytope of a degree-$5$ polynomial in two variables. Each lattice point corresponds to a monomial. The six monomials in red correspond to the vertices. The dashed triangle represents the Newton polytope of a dense polynomial of the same degree.}
    \label{fig:NewPol}
\end{figure}

Macaulay's work~\cite{macaulay,vanderwaerden} generalized the Sylvester construction to multivariate systems. This approach was later extended to \emph{sparse elimination theory}, which connects it to combinatorics, toric geometry, and hypergeometric functions~\cite{gkz1994,Sturmfels1993Sparse}. 
Sparse elimination generalizes classical results by considering the \emph{supports} of the polynomials, i.e.\ the sets of exponent vectors corresponding to monomials with nonzero coefficients.
The \emph{Newton polytope} of a polynomial, defined as the convex hull of its support (see Figure~\ref{fig:NewPol}), plays a key role in counting solutions over $(\mathbb{C}^*)^n$ via the Bernstein--Khovanskii--Kushnirenko theorem~\cite{Bernshtein1975TheNO}.
This leads to algorithms whose complexity depends on the combinatorial structure of the supports rather than total degrees.

The \emph{sparse} (or \emph{toric}) resultant generalizes the classical resultant of $n+1$ homogeneous polynomials in $n+1$ variables: it coincides with the classical resultant when all coefficients of a given degree are nonzero.
Unlike its classical counterpart, however, the sparse resultant depends only on the nonzero monomials present, typically yielding a polynomial of lower degree for sparse inputs.
Sparse elimination theory is today a standard tool in studying and solving polynomial systems, and has greatly matured since its origins in the 1970's.
This publication marks the 33th anniversary of the first paper on computational aspects of this theory, namely the one that introduced the Canny-Emiris formula~\cite{cannyemiris93}.
It is inspired by the invited talk of I.~Emiris at the Workshop on "Applications of Computer Algebra" held at Heraklion, Greece, in~2025.

This survey summarizes several results from recent decades on computing matrix representations of sparse resultants, with a special emphasis on the Canny-Emiris formula. We begin by explaining the formula and some of the tools that have recently been used  to prove the general case \cite{dandrea2020cannyemiris}. We also explain how the supports used to build the Canny-Emiris matrices are related to line bundles over toric varieties and how to build a Koszul complex whose determinant, in the sense of \cite[Appendix~A]{gkz1994}, equals the sparse resultant. Lastly, we discuss an algorithm from \cite{EmFiKoPeJ} for computing the Newton polytope of the resultant (the \emph{resultant polytope}) or its orthogonal projections, which are useful for implicitizing parametric curves and surfaces via interpolation.

\paragraph{Acknowledgments}  Checa's work is funded by the European Union under Grant Agreement number 101044561, POSALG. Views and opinions expressed are those of the authors only and do not necessarily reflect those of the European Union or European Research Council (ERC). I.~Emiris has been partially supported from the European Union’s HORIZON Europe research and innovation programme under Grant Agreement No 101136568 (HERON).

\section*{Notation}

\begin{itemize}
\item
$\kk$ is the coefficient field; unless otherwise stated,
it is arbitrary. $\PP^{n}$ denotes the projective space of dimension $n$,
obtained as a quotient of $\kk^{n+1} - \{0\}$ (non-zero multiple vectors are
identified). $\oK$ is the algebraic closure of $\kk$.
$\kk^{*}=\kk-\{0\}$ is the  field $\kk$  except $0$. 
\item $M$ will denote a lattice of rank $n$ ($\simeq \mathbb{Z}^n$), $N = \Hom_{\mathbb{Z}}(M,\mathbb{Z})$ its dual and $M_{\mathbb{R}}$ the underlying vector space. $\mathbb{T}_N = N \otimes \overline{\mathbb{K}}^*$ will denote the torus over $\overline{\mathbb{K}}$.
\item $\mathcal{A}_0,\dots,\mathcal{A}_n \subset M$ will denote a family of supports corresponding to the polynomials:
\[ 
F_i = \sum_{m \in \mathcal{A}_i}u_{i,m}\x^m \in A[x_1,\dots,x_n], \quad m \in \mathcal{A}_i \quad i = 0,\dots,n,
\]
where $\x^a = x_1^{a_1}\cdots x_n^{a_n}$ denotes a monomial (which corresponds to a  character in $\mathbb{T}_N$) associated to the lattice point $a \in \mathcal{A}_i$ and $A = \mathbb{K}[u_{i,a}, \, a \in \mathcal{A}_i, \: i = 0,\dots,n]$ denotes the ring of coefficients. Let $\Delta_i = \conv(\mathcal{A}_i) \subset M_{\mathbb{R}}$ for $i = 0,\dots,n$ be the convex hulls of the supports, also known as Newton polytopes, and $\Delta$ their Minkowski sum in $M_{ \mathbb{R}}$.

\item When we specialize $F_i$ to values of its coefficients over $\mathbb{K}$, we will denote 

\[ 
f_i = \sum_{a \in \mathcal{A}_i}\widetilde{u}_{i,a}\x^a \in \mathbb{K}[x_1,\dots,x_n] \quad a \in \mathcal{A}_i \quad i = 0,\dots,n,
\]
for $\widetilde{u}_{i,a} \in \mathbb{K}$.

\item $X_{\Sigma}$ will denote the projective toric variety associated with the normal fan  $\Sigma$ of $\Delta$. Let $R = \mathbb{K}[\X_{j}, : \rho_j \in \Sigma(1)]$ be its Cox ring and $C = R \otimes_{\mathbb{K}} A$ be the Cox ring tensored with the coefficient ring. $\F = (\F_0,\dots,\F_n)$ will denote a general homogeneous sequence of sparse polynomials in $C$ and $\f = (\f_0,\dots,\f_n)$ a  specialization. We will denote by $I(\f)$ the ideal generated by $\f$ in $R$.



\item
$\Vol_{M}(\cdot)$ is the standard Euclidean volume function in $M_{\mathbb{R}}$ and
$\MV_M(\cdot)$ denotes the mixed volume operation on arbitrary polytopes $P_1,\dots,P_n$, i.e.
$$\MV_M(P_1,\dots,P_n) = \sum_{j = 1}^n(-1)^{n-j}\sum_{0 \leq i_1 < \dots < i_j \leq n}\Vol_M(P_{i_1} + \dots P_{i_j})$$
For $i=0,1,\dots,n$, $\MV_{-i}(\Delta) = \MV(\Delta_0,\dots,\Delta_{i-1},\Delta_{i+1},\dots,\Delta_n)$ 
will denote the mixed volume
of $n$ Newton polytopes excluding the $i$-th one.
\end{itemize}

\def\0{{\mathbf{0}}}

\section{Sparse elimination}\label{sec:sparse_elim}

Elimination theory deals with the problem of finding algebraic equations for the projection of a polynomial system of equations. A typical situation is the case
of $n+1$ polynomials

\begin{equation}
\label{system}
\left\{
\begin{array}{rcl}
F_{0} & = &  \sum_{a \in \mathcal{A}_0}u_{i,a}\x^a \\
 & \vdots & \\
F_{n} & = &  \sum_{a \in \mathcal{A}_n}u_{i,a}\x^a
\end{array} 
\right.
\end{equation}

The elimination problem consists, in this case, in finding necessary and
sufficient algebraic conditions on the parameters $\widetilde{\u}= (\widetilde{u}_{i,a})_{a \in \mathcal{A}_i, i = 0,\dots,n}$ such that the specialized system of
equations $f_{0} = \ldots =  f_{n}=0$ has a root over a certain variety $X$.

Note that the solutions of the system over the torus $\mathbb{T}_N$ are invariant under translations of the $\Delta_i$. Namely, if we consider $a + \mathcal{A}_i$ as the new support of the $i$-th equation by multiplying $F_i$ by $x^a$, then $F_i(x)$ and $x^aF_i(x) = 0$ have the same zero set over $\mathbb{T}_N$. This justifies that the definition of the sparse resultant considers $\mathbb{T}_N$ as the space of solutions, as we aim to find those solutions that depend on $\Delta_0,\dots,\Delta_n$ up to translation. It might be the case that we are also interested in the solutions arising over in $\overline{\mathbb{K}}^n$ or other varieties. These solutions will arise by taking the closure of the projection of the incidence variety, or equivalently, via compactifying $\mathbb{T}_N$ using the polytopes $\Delta_0,\dots,\Delta_n$.

More precisely, consider the incidence variety $Z_F$
$$ 
\{ (\widetilde{u}_{i,a},x) \in \PP^{\mathcal{A}_0 - 1} \times \cdots \times \PP^{\mathcal{A}_n - 1}\times \mathbb{T}_N :
\: f_0 = \dots = f_n = 0\}.
$$
Here, $\mathbb{P}^{\mathcal{A}_i - 1}$, the projective space of dimension $|\mathcal{A}_i| - 1$, is the space of coefficients of the $i$-th equation, whose zero set is invariant under scalar multiples.

The incidence variety has a natural projection
\begin{eqnarray*}
\pi : Z_F & \rightarrow & 
\PP^{\mathcal{A}_0 - 1} \times \cdots \times \PP^{\mathcal{A}_n - 1}.
\end{eqnarray*}
The Zariski closure of $\pi(Z_F)$ equals the set of parameters $(\widetilde{u}_{i,a})_{a \in \mathcal{A}_i, i = 0,\dots,n}$ for
which the system \eqref{system} can be homogenized so that it has a root over the compactification $X_{\Sigma}$ determined by the polytopes $\Delta_0,\dots,\Delta_n$. This compactification equals a toric variety, as we shall see in \Cref{section:toricvarieties}. The goal of sparse elimination is to find algebraic equations for $\pi(Z_F)$.

\begin{definition}
\label{sparse:resultant}
The \textit{sparse eliminant}, denoted $\Elim_{\mathcal{A}}$ is the irreducible polynomial defining the Zariski closure of  $\pi(Z_F)$, if $\codim(\pi(Z_F)) = 1$, or $1$ otherwise. The \textit{sparse resultant}, denoted $\Res_{\mathcal{A}}$ is any primitive equation defining the push-forward $\pi_*(Z_F)$.
\end{definition}

\begin{example}
If we consider the classical case where $\mathcal{A}_i$ consists of all monomials of degree at most $d_i$, for a sequence $d_0 \leq \dots \leq d_n$, then the zeros of the resultant correspond to the coefficients of homogeneous polynomial systems that have a root over $\PP^{n}$.
\end{example}

The resultant is always a power of the eliminant, i.e.,
\begin{equation}
\label{resultanteliminant}
\Res_{\mathcal{A}} = \Elim_{\mathcal{A}}^{d_{\mathcal{A}}}, \quad d_{\mathcal{A}} \in \mathbb{Z}_{>0}.
\end{equation}
 
 As they often coincide, the resultant and the eliminant are sometimes cross-refered in the literature  \cite{sturmfels94}. However, the distinction between these two definitions, appearing in \cite{esterov2010newton, poissondansom}, clarifies which object we are computing: while for most practical purposes, computing the eliminant is sufficient, the object computed by some of the closed formulas we discuss below is the sparse resultant. 

Definition \ref{sparse:resultant} can be given in other settings where the $F_i$ are not necessarily given by a monomial structure, but could depend on other families of polynomial functions. Mutatis mutandis, a definition of the eliminant can be given in a similar manner over an irreducible projective variety, instead of $X_{\Sigma}$  \cite{BUSE2000515, MONIN2020107147,Kapranov1992ChowPA}.

\begin{definition}\label{dfn:essential}
Let $J \subset \{0,\dots,n\}$ be a subset and consider the lattice $\mathcal{L}_j$ spanned by the differences of lattice points in $\mathcal{A}_j$. Let $\mathcal{L}_J = \sum_{j \in J}\mathcal{L}_j$ and $\rank(J) = \rank(\mathcal{L}_J)$. A subset of supports $(\mathcal{A}_i)_{i \in I}$ is essential if 
$$\rank(I) = |I| -  1, \quad \rank(J) \geq |J|, \: \forall J \subsetneq I.$$
\end{definition}

The following theorem, due to Sturmfels \cite[Theorem 1.1, Corollary 1.1]{sturmfels94}, characterizes the case where $\pi(Z_F)$ has codimension $1$.

\begin{theorem}
\label{codimensionone}
    The codimension of $\pi(Z_F)$ equals the maximum of $|J| - \rank(J)$ over all $J \subset \{0,\dots,n\}$. In particular, $\codim(\pi(Z_F)) = 1$ if and only if there exists a unique essential family.
\end{theorem}

Note that the exponent $d_{\mathcal{A}}$ in \eqref{resultanteliminant} can be expressed by saturating the lattice $L_J$ where $J$ is the unique essential family \cite[Proposition 3.6]{dandrea2020cannyemiris}. 
In the cases where the essential subfamily is strictly smaller than $\{0,\dots,n\}$, the resultant  depends only on the coefficients of those polynomial equations and the overdetermined system lies over a lattice of strictly smaller dimension.

\begin{example}
    Consider the following overdetermined system
    $$F_0 = u_{00} + u_{01}x, \: F_1 = u_{10} + u_{11}y, \: F_2 = u_{20} + u_{21}y.$$
    Note that $J = \{1,2\}$ is the unique essential family. Thus, the resultant  depends only on the coefficients of $F_1,F_2$. In this case, it is straightforward to check that $\Res_{\mathcal{A}} = \Elim_{\mathcal{A}} = u_{10}u_{21} - u_{11}u_{20}$.
\end{example}

Note that the problem of sparse elimination, can be simplified under some assumptions on the supports $\mathcal{A}_i$. 
\begin{itemize}
    \item When all the supports coincide, we say that the system is unmixed. Computing the sparse resultants for unmixed systems with support $\mathcal{A}$ is equivalent to computing the Chow form of the toric variety defined by the monomial parametrization of $\mathcal{A}$ \cite{dalbecsturmfels,sturmfelspedersen}. 
    When the number of supports is very close to the number of variables or, equivalently, when the embedded toric variety defined by $\mathcal{A}$ has low codimension, formulas for Chow forms can be made quite explicit \cite{DICKENSTEIN2002119, Sturmfels1993Sparse}.  Resultants of bivariate unmixed systems have also been studied in \cite{KHETAN2003425, KHETAN2005237, CHTCHERBA2004915}.
    \item Other types of structured supports, such as star multilinear systems or generalized unmixed systems, have also been considered \cite{dandrea2002, benderthesis}.
    \item For the dense resultant, the original formula of Macaulay \cite{macaulay}, which relied on the monomials of a certain degree, has been shown to be a particular case of the subdivision-based formula that we discuss here \cite[Section 5]{dandrea2020cannyemiris}. Other similar formulas, involving Bezout resultants or Morley forms, appear in \cite{Jouanolou1997FormesDE} and approaches for a sparse generalization have appeared in \cite{BUSE2024107739, Eisenbud2003Resultants, weymannzelevinsky}.
    \item As it is the first step in generalizing the dense case, formulas for the multihomogeneous resultant, where the supports are dense in the monomials of a fixed multidegree in several groups of variables, have also received wide attention \cite{dickensteinemirismultihomo03, emimantzdeterminant, bender2021koszultype, BUSE2022514, BUSE202065}. 
\end{itemize}

\section{The Canny-Emiris formula}
\label{section:cannyEmiris}
In this section, we recall a formula for computing $\ResA$ which was introduced in \cite{cannyemiris93} and proved (under some assumptions) in \cite{dandrea2020cannyemiris}. This formula relies on the construction of a mixed subdivision of $\Delta$.

\begin{definition}
\label{infconvolution}
A mixed subdivision of $\Delta$ is a decomposition of this polytope in a union of cells $\D$ such that $\Delta = \cup \D$ and:
\begin{itemize}
    \item[i)] the intersection of two cells is either a cell or empty,
    \item[ii)] every face of a cell is also a cell of the subdivision and,
    \item[iii)]  every cell $\D$ has a component structure $\D = \D_0 + \dots + \D_n$ where $\D_i$ is a cell of a subdivision in $\Delta_i$.
    Moreover, these component structures of the cells are compatible, i.e.\ coincide on a common face of a pair of cells.
\end{itemize}
\end{definition}

The usual way to construct mixed subdivisions is to consider piecewise affine convex lifting functions $\rho_i: \Delta_i \xrightarrow[]{} \mathbb{R}$ and to take the lower convex hull of the Minkowski sum of the graphs of each of the functions. A piecewise affine global lifting function $\rho: \Delta \xrightarrow[]{} \mathbb{R}$ is obtained after taking the inf-convolution of the previous functions, i.e.
$$\rho(x) = \inf \{\rho_0(x_0) + \dots \rho_n(x_n) \: : x = x_0 + \dots + x_n\}.$$

\begin{figure}
\vspace{-0.9cm}

$$   \includegraphics[scale=0.35]{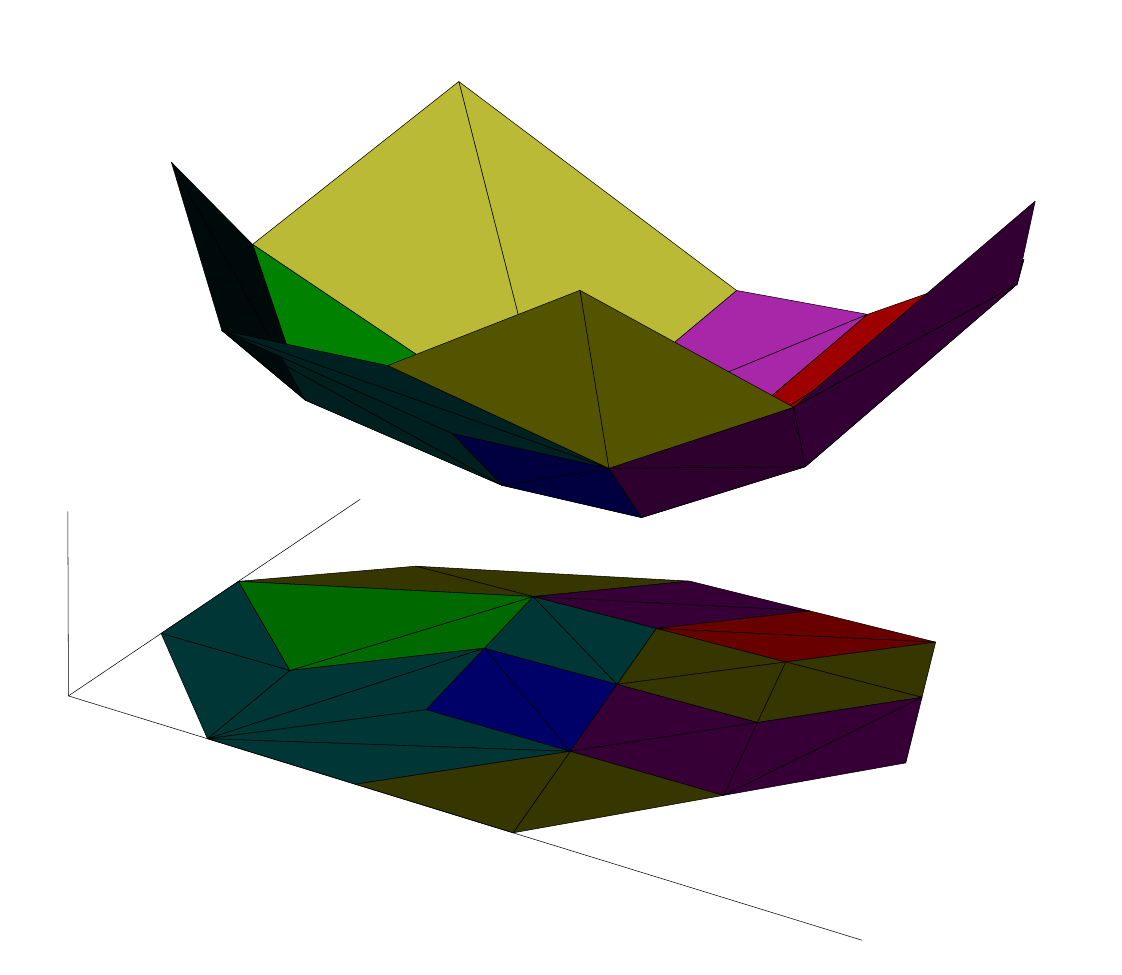}\includegraphics[scale=0.35]{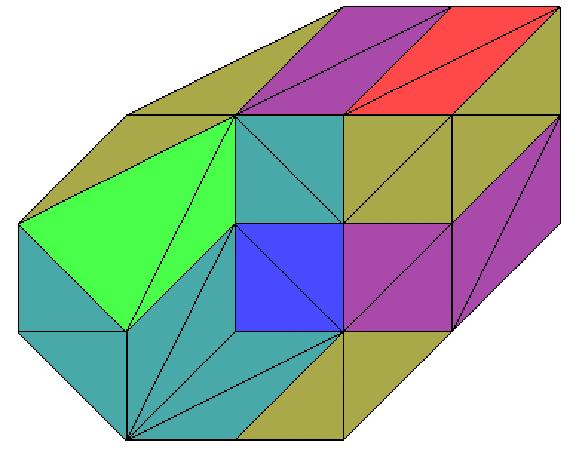}$$
  \caption{The way to construct coherent mixed subdivisions is considering piecewise affine convex lifting functions $\rho_i: \Delta_i \xrightarrow[]{} \mathbb{R}$. Then, take the Minkowski sum of the graphs of these functions as polytopes in $\mathbb{R}^{n+1}$ and project the lower convex hull back to $\mathbb{R}^n$. Each of the faces of the lower convex hull defines a cell in the subdivision.}
  \label{fig:mixedsub}
\end{figure}

We denote the mixed subdivision constructed this way as $S(\rho)$. The cells $\D \in S(\rho)$ corresponding to the subsets of $\Delta$ where $\rho$ is affine correspond to the cells of a mixed subdivision. In order to define $\rho_i$, it suffices to define the values of this function on the lattice points $a \in \mathcal{A}_i$ which we denote as $\omega_{i,a}$. We identify the lifting function $\rho$ with its values $\omega \in \prod_{i = 0}^n\mathbb{R}^{\mathcal{A}_i}$. Subdivisions of this type are called \textit{coherent} (or \textit{regular}) (see \Cref{fig:mixedsub}).

This technique of deformation is commonly also used for square systems. For example, it can be used to describe parameter homotopies \cite{hubersturmfels}, study real algebraic varieties \cite{stuviro}, find monomial bases of the quotient ring \cite{mixedmonomial, emirisrege}, decide upon the nonnegativity of polynomials \cite{feliu2025copositivitydiscriminantsnonseparablesigned}, or count the number of connected components of fewnomial surfaces \cite{bihan2024newboundsnumberconnected}. 
In all these cases, one considers polynomials of the form

\begin{equation}
\label{ft}
F_t = \sum_{m \in \mathcal{A}}u_{m}\x^{m}t^{\omega_m},
\end{equation}
for some finite set of supports $\mathcal{A} \subset M$ and $\omega_{m} \in \mathbb{R}$.

\begin{definition}
A mixed subdivision of $\Delta$ is tight if, for every $n$-dimensional cell $F$, its components satisfy:
\[ \sum_{i=0}^n\dim \mathcal{D}_i = n.\]
\end{definition}
In the case of $n + 1$ polynomials and $n$ variables, this property guarantees that every $n$-dimensional cell has a component that is $0$-dimensional. Cells that have a single $0$-dimensional component are called \textit{mixed} ($i$-mixed if it is the $i$-th component). The rest of the cells are called \textit{non-mixed}. Let $\delta$ be a generic vector such that the lattice points in the interior of $\Delta + \delta$ lie in $n$-dimensional cells. Then, consider:
 \begin{equation}
 \label{banddelta}
 \mathcal{B} = (\Delta + \delta) \cap M.\end{equation}
Each element $b \in \mathcal{B}$ lies in one of the translated cells $\D + \delta$ and let $\D_i$ be the components of this cell. As the subdivision is tight, there is at least one $i$ such that $\dim \D_i = 0$.

 \begin{definition}
 \label{rowcontent}
 The \textit{row content} of $S(\rho)$ is a function 
 \[ 
 \rc: \mathcal{B} \xrightarrow[]{} \cup_{i = 0}^n\{i\} \times \mathcal{A}_i \] where, for $b \in \mathcal{B}$ lying in an  $n$-dimensional cell $\D$, $\rc(b)$ is a pair $(i(b),a(b)) \in \cup_{i = 0}^n \{i\} \times \mathcal{A}_i$ with $$i(b) = \max\{i \in \{0,\dots,n\} \: | \: b \in \D, \: \dim \D_{i} = 0\}, \quad a(b) = F_{i(b)}.$$
 \end{definition}
This provides a partition of $\mathcal{B}$ into subsets:
 \[ \mathcal{B}_i = \{b \in \mathcal{B} \: : \: i(b) = i\}.
 \]
The row content function is defined in such a way that $b - a(b) + \mathcal{A}_{i(b)} \subset \mathcal{B}$ for all $b \in \mathcal{B}$ \cite[Lemma 5.4]{cannyemiris}.

With a row content function, one can construct Canny-Emiris matrices $\mathcal{H}_{\mathcal{A},\rho}$ whose rows correspond to the coefficients of the polynomials $\x^{b - a(b)}F_{i(b)}$ for each of the $b \in \mathcal{B}$. The rows and columns of $\mathcal{H}_{\mathcal{A},\rho}$ are labeled by $\mathcal{B}$ and the entry corresponding to a pair $b,b' \in \mathcal{B}$ equals
\[\mathcal{H}_{\mathcal{A},\rho}[b,b'] =  \begin{cases} 
      u_{i(b),b'-b+a(b)} & b' - b + a(b) \in \mathcal{A}_i \\
      0 & \text{otherwise}
   \end{cases}
\]

Each entry contains, at most, a single coefficient $u_{i,a}$. In particular, the row content allows us to choose a maximal submatrix of $\mathcal{H}_{\mathcal{A},\rho}$ of the matrix of the \textit{Macaulay map} $\mathcal{M}_{\mathcal{A}}$
\begin{multline}
\label{macaulay}
\bigoplus_{i = 0}^n V_i \xlongrightarrow{\mathcal{M}_{\mathcal{A}}} V, \quad (G_0,\dots,G_n) \xrightarrow[]{} \sum_{i = 0}^nG_iF_i, \quad  V_i = \bigoplus_{b \in \mathcal{B}_{(i)}} \mathbb{K}\cdot \x^{b}, \quad   V = \bigoplus_{b \in \mathcal{B}} \mathbb{K}\cdot \x^{b} \end{multline}
where $\mathcal{B}_{(i)} = \{b - \mathcal{D}_i : \: b \in \mathcal{B} \cap \mathcal{D}\}$. This class of matrices are called Sylvester-type matrices, as they generalize the construction of Sylvester for the resultant of two univariate polynomials, and contain a single coefficient in each entry.

Let $\mathcal{C} \subset \mathcal{B}$ be a subset of the supports in translated cells. A submatrix $\mathcal{H}_{\mathcal{A},\rho,\mathcal{C}}$ is defined by considering the submatrix of the corresponding rows and columns associated with elements in $\mathcal{C}$. In particular, the set of lattice points lying in translated non-mixed cells is defined as
\[ \mathcal{B}^{\circ} \coloneqq \{b \in \mathcal{B} \: | \:  b \text{ lies in a translated non-mixed cell}\}.\]
$\mathcal{B}^{\circ}$ indexes the principal submatrix $\mathcal{E}_{\mathcal{A},\rho} = \mathcal{H}_{\mathcal{A},\rho,\mathcal{B}^{\circ}}$.

The Canny-Emiris conjecture states that the sparse resultant is the quotient of the determinants of these two matrices:
\begin{equation}
\label{rationalformula}
\Res_{\mathcal{A}} = \frac{\det(\mathcal{H}_{\mathcal{A},\rho})}{\det(\mathcal{E}_{\mathcal{A},\rho})}.\end{equation}

\begin{remark}
    One can also construct bigger or smaller Canny-Emiris matrices than $\mathcal{H}_{\mathcal{A},\rho}$ with the following two strategies:
    \begin{itemize}
        \item[$i)$] Consider $\Lambda \subset M_{\mathbb{R}}$ a polytope and replace $\mathcal{B}$ by $(\Delta + \Lambda) \cap M$. A tight mixed subdivision of $\Delta + \Lambda$ also induces a row content function as in Definition \ref{rowcontent}. Although this strategy will produce bigger matrices, which is not convenient for computational purposes, their existence will be helpful for the explanations in Section \ref{section:toricvarieties}.
        
        \item[$ii)$] If there exists a subset $\mathcal{G} \subset \mathcal{B}$ such that $b - a(b) + \mathcal{A}_{i(b)} \subset \mathcal{G}$ and $\mathcal{G}$ contains all lattice points in mixed cells, then the submatrix $\mathcal{H}_{\mathcal{A},\rho,\mathcal{G}}$ that considers only the lattice points in $\mathcal{G}$ can be used. This leads to the greedy algorithm in \cite{cannypedersen}.
    \end{itemize}
\end{remark}

\begin{example}
\label{smallexample}
Let $F_0,F_1,F_2$ be three bilinear equations corresponding to the supports $\mathcal{A}_0 =  \mathcal{A}_1 =  \mathcal{A}_2 = \{(0,0),(1,0),(0,1),(1,1)\}$. A possible mixed subdivision $S(\rho)$ is the following:

\[\begin{tikzpicture}[scale=1.7, transform shape]

\draw[brown] (1/2,1/2) -- (1,1/2);
\draw[brown] (1/2,1/2) -- (1/2,1);
\draw[brown] (1,1/2) -- (1,1);
\draw[brown] (1/2,1) -- (1,1);
\draw[brown] (1/2,0) -- (1,0);
\draw[brown] (0,1/2) -- (0,1);
\draw[brown] (1/2,3/2) -- (1,3/2);
\draw[brown] (3/2,1/2) -- (3/2,1);

\draw[blue] (0,0) -- (1/2,0);
\draw[blue] (0,0) -- (0,1/2);
\draw[blue] (1/2,0) -- (1/2,1/2);
\draw[blue] (0,1/2) -- (1/2,1/2);
\draw[blue] (1,1/2) -- (1,0);
\draw[blue] (1/2,1) -- (0,1);
\draw[blue] (3/2,1/2) -- (3/2,0);
\draw[blue] (1/2,3/2) -- (0,3/2);

\draw[green] (1,1) -- (3/2,1);
\draw[green] (1,1) -- (1,3/2);
\draw[green] (3/2,3/2) -- (3/2,1);
\draw[green] (3/2,3/2) -- (1,3/2);
\draw[green] (3/2,0) -- (1,0);
\draw[green] (0,3/2) -- (0,1);
\draw[green] (3/2,1/2) -- (1,1/2);
\draw[green] (1/2,3/2) -- (1/2,1);

\filldraw[black] (0+1/4,1/4) circle (2pt) node[anchor=east] {};
\filldraw[red] (0+1/4,1/2+1/4) circle (2pt) node[anchor=east] {};
\filldraw[red] (0+1/4,1+1/4) circle (2pt) node[anchor=east] {};
\filldraw[red] (1/2+1/4,0+1/4) circle (2pt) node[anchor=south ] {};
\filldraw[red] (1/2+1/4,1/2+1/4) circle (2pt) node[anchor= south] {};
\filldraw[red] (1/2+1/4,1+1/4) circle (2pt) node[anchor=west] {};
\filldraw[red] (1+1/4,0+1/4) circle (2pt) node[anchor=south] {};
\filldraw[red] (1+1/4,1/2+1/4) circle (2pt) node[anchor=west] {};
\filldraw[red] (1+1/4,1+1/4) circle (2pt) node[anchor=west] {};

\end{tikzpicture}\]
where the dots indicate the set $\mathcal{B}$ of lattice points  in translated mixed cells, and the colors indicate the support $\Ac_i$ to which the summand belongs.  The number of lattice points in translated cells is $9$. However, if we construct the matrix greedily starting from the lattice points in translated mixed cells, we have an $8 \times 8$ matrix, corresponding to the lattice points marked in red.
\end{example}

The proof of the formula in \eqref{rationalformula} given by D'Andrea, Jerónimo, and Sombra in \cite{dandrea2020cannyemiris} relies on imposing that the mixed subdivision $S(\rho)$ given by the lifting $\rho$ satisfies a certain condition, given on a chain of mixed subdivisions.

\begin{definition}
Let $S(\phi), S(\psi)$ be two mixed subdivisions of $\Delta = \sum_{n = 0}^n\Delta_i$. We say that $S(\psi)$ refines $S(\phi)$ and write $S(\phi) \preceq S(\psi)$ if for every cell $\mathcal{D}' \in S(\psi)$ there is a cell $\D \in S(\phi)$ such that $\mathcal{D}' \subset \D$ and  $\mathcal{D}_i' \subset \D_i$, for all $i$. An incremental chain of mixed subdivisions $S(\theta_0) \preceq \dots \preceq S(\theta_n)$ is a chain of mixed subdivisions of $\Delta$ ordered by refinement. 

\end{definition}

The tropical analogue of refinement of mixed subdivisions was studied in \cite{mcs}.

\begin{remark}
\label{remarkzerolifting}
In \cite[Definition 2.4]{dandrea2020cannyemiris}, a common lifting function $\omega \in \prod_{i = 0}^n\mathbb{R}^{\mathcal{A}_i}$ is considered and the $S(\theta_i)$ are given by the lifting functions 
\[
\omega^{<i} = (\omega_0,\dots,\omega_{i-1},\0)
\] 
as long as $S(\theta_i) \preceq S(\theta_{i+1})$. The last zero $\0$ represents the lifting on $(\Delta_i,\dots,\Delta_n)$ for $i = 0,\dots,n$. The resulting mixed subdivision is the same as if we considered the zero lifting in $\sum_{j = i}^n\Delta_j$.
\end{remark}

The following theorem can be found in \cite[Theorem 3.4]{emirisrege}

\begin{proposition} 
\label{mixedvolume1}
Let $S(\rho)$ be a tight mixed subdivision of $\Delta = (\Delta_0,\dots,\Delta_n)$. For $i = 0,\dots,n$, the mixed volume of all polytopes except $\Delta_i$ equals the volume of the $i$-mixed cells.
\[\MV_{-i}(\Delta) = \sum_{\mathcal{D} \: i\text{-mixed}}\Vol_n \mathcal{D}\]
\end{proposition}
In particular, $\MV_{-i}(\Delta)$ equals the degree of the sparse resultant in the coefficients of $F_i$ \cite[Chapter 7, Theorem 6.3]{coxlitosh}. Each of the rows of $\mathcal{H}_{\mathcal{A},\rho}$ will correspond to a lattice point $b$, and each entry in that row will have degree $1$ with respect to the coefficients of $F_{i(b)}$, and zero with respect to the coefficients of the rest of polynomials. Therefore, if we add the lattice points in $i$-mixed cells, the degree of $\mathcal{H}_{\mathcal{A},\rho}$ with respect to the coefficients of $F_i$ will be at least the degree of the resultant with respect to the same coefficients.

From Proposition \ref{mixedvolume1}, it follows that if the unique essential family is $J = \emptyset$, then the resultant is equal to $1$. Moreover, if the unique essential family equals $\{i\}$ then $\mathcal{A}_i$ is given by a single point $\{a\}$ and the sparse resultant  equals $u_{i,a}^{m_i}$ for $m_i = \MV(\Delta_0,\dots,\Delta_{i-1},\Delta_{i+1},\dots,\Delta_n)$. The Canny-Emiris formula holds \cite[Proposition 4.26]{dandrea2020cannyemiris} in both cases. In the following proposition, $\mathcal{B}_{\mathcal{D},i} \coloneqq \{b \in \mathcal{B} \cap (\mathcal{D}+\delta) : \: i(b) = i\}$. 
\begin{definition}
\label{admissible}
An incremental chain $S(\theta_0) \preceq \dots \preceq S(\theta_n)$ is admissible if for each $i = 0,\dots,n$, each $n$-dimensional cell $\D$ of the subdivision $S(\theta_i)$ satisfies either of the following two conditions
\begin{itemize}
    \item[$i)$] the unique essential family of $\mathcal{A}_{\D}$ contains at most one support,
    \item[$ii)$] $\mathcal{B}_{\D,i}$ is contained in the union of the translated $i$-mixed cells of $S(\rho_{\D})$.
\end{itemize}
A mixed subdivision $S(\rho)$ is called \textit{admissible} if it admits an admissible incremental chain $S(\theta_0) \preceq \dots \preceq S(\theta_n) \preceq S(\rho)$ refining it.
\end{definition}

While in case $(i)$, the mixed subdivision the formula holds for the supports $\mathcal{A}_{\mathcal{D}} = (\mathcal{A}_0 \cap \mathcal{D},\dots,\mathcal{A}_n \cap \mathcal{D})$, 
case $(ii)$ is useful to prove the formula by induction on the incremental chain. 
The last ingredients of this proof are the following two product formulas for the initial parts of $\Res_{\mathcal{A}}$ and $\det(\mathcal{H}_{\mathcal{A},\rho})$, which appear in \cite{dandrea2020cannyemiris}.

\begin{theorem}
\label{initialFormTheorem}
Let $S(\phi) \preceq S(\rho)$ be a chain of mixed subdivisions. For any polynomial $F_t$ as in \eqref{ft}, let $\init_{\omega}(F_t)$ be the (sum of) coefficients of the terms of $F_t$ of minimal degree in $t$. Then,
\begin{equation}
        \label{equation1}
\init_{\omega}(\Res_{\mathcal{A}}(u_{i,a}t^{\omega_{i,a}})) = \prod_{\mathcal{D} \in S(\rho)} \Res_{\mathcal{A}_{\mathcal{D}}}
\end{equation}\
    where $\mathcal{A}_{\D}$ refers to the set of supports $(\mathcal{A}_0 \cap \D_0,\dots,\mathcal{A}_n \cap \D_n)$. Moreover, if $\widetilde{\omega} \in \prod_{i = 0}^n\mathbb{R}^{\mathcal{A}_i}$ are the values defining $S(\phi)$, then
    \begin{equation}
    \label{equation2}
    \init_{\widetilde{\omega}}(\det(\mathcal{H}_{\mathcal{A},\mathcal{\rho}}(u_{i,a}t^{\widetilde{\omega}_{i,a}})) = \prod_{\mathcal{D} \in S(\phi)}\det(\mathcal{H}_{\mathcal{A}_D,\mathcal{\rho}_D}). \end{equation}
\end{theorem}

For generic $\omega$, the polynomial in \eqref{equation1} is a monomial corresponding to a vertex of $N(\ResA)$. In \Cref{sectionpolytopes}, we recall a formula for this monomial, which will be useful in computing $N(\ResA)$. The proof of \eqref{equation2} follows by noting that the minimum of the function
\begin{equation}
\label{equation_cannyemiris}
\rho(b-\delta) - \rho_{i(b)}(a(b)) + \rho_{i(b)}(a'), \quad a' = b' - b + a(b), \quad b \in \mathcal{B}, b' \in b - a(b) + \mathcal{A}_{i(b)} \end{equation}
is only attained at $b = b'$ \cite[Lemma 4.11]{dandrea2020cannyemiris}. This allowed the authors to generalize an idea appearing also in \cite{sturmfels94, cannyemiris}: the terms in the diagonal of matrix $\mathcal{H}(u_{i,m}t^{\omega_{i,m}})$ are minimal with respect to the degree in $t$, and so the initial part is the product of the diagonal terms. With all of the above ingredients, a proof of the rational formula \eqref{rationalformula} was given, provided that $S(\rho)$ admits an admissible incremental chain refining it; see \cite[Theorem 4.27]{dandrea2020cannyemiris}. An approach providing combinatorial conditions enforcing the existence of an admissible chain refining $S(\rho)$ was given in \cite{mcs}.

\section{Toric varieties, Koszul complexes and resultants}
\label{section:toricvarieties}

In this section, we explain the toric geometry counterpart of the formula described in \eqref{rationalformula}. The treatment of some of the aspects in this chapter can be found in \cite{bender2021toric, sparsenullstellensatz, BUSE2024107739}, although the relation between resultants, complexes and toric varieties goes back to \cite{gkz1994}. A much more general treatement appears in the last chapter of  \cite{weymancohomo}, where it is suggested that all possible formulas for resultants arise as determinants of a general complex.

Let us consider the normal fan $\Sigma$ of the Minkowski sum $\sum_{i = 0}^n \Delta_i$, see \Cref{sec:intro}. The underlying complete toric variety $X_{\Sigma}$ can be defined as the closure of the image of a monomial map, or directly using the fan $\Sigma$ \cite{coxlittleschneck}.

For square systems, the toric variety $X_\Sigma$ provides a natural geometric setting as the root count over $\mathbb{T}_N$ coincides with the intersection number of $\mathbb{T}_N$-invariant divisors on this variety, which are associated to their Newton polytopes. Moreover, the toric vector bundle defined as the direct sum of line bundles can be attached to a Koszul complex. We can recover the sparse resultant as the determinant of this complex at some degrees.

Let $\Sigma(1) = \{\rho_j\}_{j = 1,\dots,n+r}$ be the set of rays of $\Sigma$ (the $1$-dimensional cones) and let $D_j$  be the corresponding $\mathbb{T}_N$-invariant divisor corresponding to that ray, via the orbit-cone correspondence \cite[Theorem 3.2.6]{coxlittleschneck}. 

Each $\mathbb{T}_N$-invariant divisor is of the form $D = \sum_{j = 1}^{n + r}a_jD_j$ and corresponds to a class $\alpha = [D] \in \Cl(X_{\Sigma})$. Each of these divisors can be associated to a polytope
\begin{equation}\label{eq:polytopes}
    \Delta_{D} = \bigl\{ m \in M_{\mathbb{R}} \; : \; \langle u_j, m \rangle \geq -a_{j},\; j = 1,\dots,n+r \bigr\} \qquad (a_j) \in \mathbb{Z}^{n+r}
\end{equation}
where $u_j$ is the primitive generator of the ray $\rho_j$. In particular, each polytope whose normal fan is refined by $\Sigma$ has an associated divisor $D$ constructed from its face description. Although for the following construction it is not critical to understand the main properties of $X_{\Sigma}$, it is crucial to understand the definition of the Cox ring.

\begin{figure}
    \centering
\medskip
    
    \includegraphics[width=0.5\linewidth]{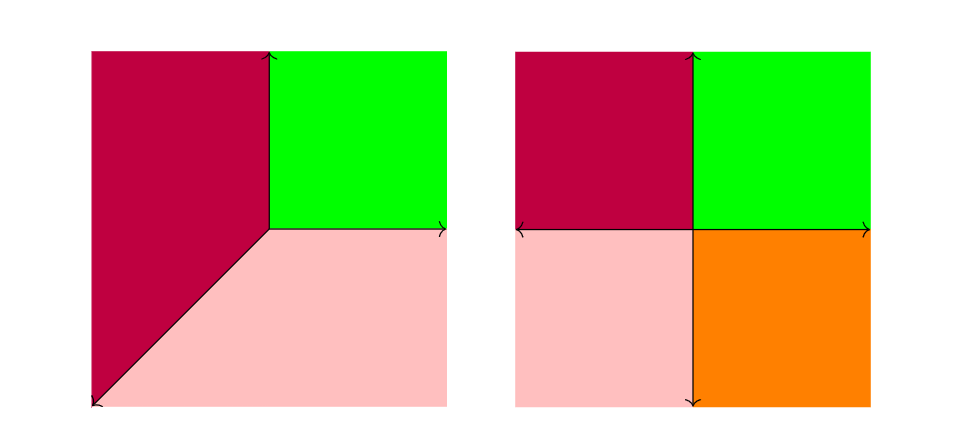} 
    \caption{The normal fans of a triangle and a square in $\mathbb{R}^2$. They define the complete toric varieties $\PP^2$ and $\PP^1 \times \PP^1$, respectively.}
    \label{fig:placeholder}
\end{figure}

\begin{definition}
    The \textit{Cox ring} (or homogeneous coordinate ring) of $X_{\Sigma}$ is the ring
    $$R \coloneqq \mathbb{K}[\mathbf{X}_{j} \: : j = 1,\dots,n + r]$$
with the $\Cl(X_{\Sigma})$-grading in $R$ is given by $\deg(\X_j) = \pi(e_j)$, where $\pi$ is the second map of the short exact sequence
    \begin{equation}
    \label{ses}
        0 \xrightarrow[]{} M \xrightarrow[]{\mathbf{U}} \mathbb{Z}^{n + r} \xrightarrow[]{\pi} \Cl(X_{\Sigma}) \xrightarrow[]{} 0
    \end{equation}
    where $\mathbf{U}$ is the matrix given by $(u_j)_{j = 1,\dots,n+r}$ which are the primitive generators of the rays $\rho_j \in \Sigma(1)$.
\end{definition}

 In particular,
\begin{equation}
\label{monomialbasis}
R_{\alpha} \simeq \bigoplus_{m \in \Delta_{D}} \mathbb{K}\cdot x^m, \end{equation}
for any $\alpha \in \Cl(X_{\Sigma})$ and $[D] = \alpha$. We denote by $R(-\alpha)$ the ring $R$ with the grading given by $R(-\alpha)_{\alpha'} = R_{\alpha'-\alpha}$ for $\alpha,\alpha' \in \Cl(X_{\Sigma})$.

\begin{remark}
\label{codimension}

Under the condition of codimension~1 in Theorem~\ref{codimensionone}, we can always assume that $\Delta$ is $n$-dimensional. Otherwise, we can restrict to the unique essential family and work in the lattice spanned by those supports.
\end{remark}

Following \cite[§6.1]{coxlitosh}, note that for $D = \sum_{j = 1}^{n+r}a_jD_j$, there exist $m_{\sigma}$ such that
$$\langle u_j, m_{\sigma}\rangle = -a_j, \quad \rho_j \in \sigma.$$
The definitions of nef and ample divisors, which come from algebraic geometry, have a very direct combinatorial interpretation in terms of the polytopes \cite[Theorem 6.1.7, Theorem 6.1.13]{coxlitosh}.

\begin{theorem}
\label{nefample}
    Let $D = \sum_{j = 1}^{n +r}a_jD_j$ be a $\mathbb{T}_N$-invariant divisor for
    \begin{itemize}
        \item $D$ is nef, if and only if, $m_{\sigma} \in \Delta_{D}$ for all $\sigma \in \Sigma$. 
        \item $D$ is ample, if and only, $\langle u_j, m_{\sigma} \rangle > -a_{j}$ for all $\rho_j \not\subset \sigma$.
    \end{itemize}
\end{theorem}

We denote by $\Nef(X_{\Sigma})$ the cone of nef Cartier divisors.

\begin{remark}
    The last item in Theorem \ref{nefample} is also equivalent to the fact that the normal fan of $\Delta_{D}$ equals $\Sigma$.
\end{remark}

Next, we return to the setting of $n + 1$ polytopes in $M_{\RR}$. Let $D_0,\dots,D_n$ be nef Cartier divisors, defined as
$$D_i = \sum_{j = 1}^{n + r}a_{ij}D_j, \qquad a_i \coloneqq (a_{ij}) \in \mathbb{Z}^{n +r}$$
associated with polytopes $\Delta_0,\dots,\Delta_n$. A \textit{general homogeneous polynomial system} $\F = (\F_0,\dots,\F_n)$ is the system defined by the polynomials 
\begin{equation}\label{eq:Fi}
\F_i = \sum_{\X^{\mu} \in R_{\alpha_i}}c_{i,\mu}\mathbf{X}^\mu \in C, \ \ i=0,\ldots,n.
\end{equation}
Note that $\F_i$ corresponds to a general global section of $\mathcal{O}(D_i)$ and so $\F = (\F_0,\dots,\F_n)$ is  a general global section of $\oplus_{i = 0}^n \mathcal{O}(D_i)$.

Choosing an $n$-dimensional cone $\sigma \in \Sigma$, one can dehomogenize this system by setting the variables $\X_j$ for $\rho_j \not\subset \sigma$ to $1$. If the cone coincides with the first orthant in $M_{\mathbb{R}}$, one recovers a polynomial system with supports equal to $\mathcal{A}_i = \Delta_i \cap M$. Conversely, the polynomials
 $\F_0,\dots,\F_n$ can be defined as the homogenizations of the polynomials  $F_0, \ldots , F_n$ with supports in the subsets $\mathcal{A}_i = \Delta_i \cap M$, $i = 0,\dots,n$. More precisely, the homogenization of the polynomial 
\begin{equation}
    \label{nothomogenized}F_i = \sum_{m \in \mathcal{A}_i}u_{i,m}x^m \in  A[x_1,\dots,x_n]
\end{equation}
is the polynomial 
\begin{equation}
    \label{eqhomogenized}\F_i = \sum_{m \in \mathcal{A}_i}u_{i,m}\X^{\mathbf{U}m + a_i} \in {C} = A[\X_1,\dots,\X_{n + r}]
\end{equation}
where $\mathbf{U}$ is defined in \eqref{ses}.

\begin{definition}
Let $E = \langle e_0,\dots,e_n\rangle $ be a vector space of dimension $n + 1$. The \textit{Koszul complex} associated with $\F = (\F_0,\dots,\F_n)$ is the complex of $R$-modules
    \begin{equation}\label{koszulcomplex}
K_{\bullet}(\F): \bigoplus_{|I| = n + 1}R(-\alpha_I) \otimes e_I \xrightarrow[]{\partial_{n+1}} \dots \xrightarrow[]{\partial_3} \bigoplus_{|I| = 2}R(-\alpha_I)\otimes e_I \xrightarrow[]{\partial_2} \bigoplus_{i = 0}^nR(-\alpha_i)\otimes e_i  \xrightarrow{\partial_1} R.
\end{equation} for $\alpha_I = \sum_{i \in I}\alpha_i$ and $e_I = \wedge_{i \in I}e_i \in \bigwedge^{|I|}E$. The differentials are defined as
$$\mathbf{\partial}_i(\G \otimes e_I) = \sum_{i \in I}(-1)^{\#(i)}\G \F_i \otimes e_{I \setminus \{i\}}$$
where $\#(i)$ denotes the position of $i$ in the ordered set $I$, and extend linearly. We denote $K_{\bullet}(\f)$ the Koszul complex when $\f_i$ are specializations of the $\F_i$. 
\end{definition}

This complex of $R$-modules is graded and its graded pieces are the complexes of 

\begin{equation}\label{koszulcomplexgraded}
K_{\bullet}(\F)_{\alpha}: \bigoplus_{|I| = n + 1}R_{\alpha - \alpha_I} \xrightarrow[]{\partial_{n+1}} \dots \xrightarrow[]{\partial_3} \bigoplus_{|I| = 2}R_{\alpha - \alpha_I} \xrightarrow[]{\partial_2} \bigoplus_{i = 0}^nR_{\alpha - \alpha_I}  \xrightarrow{\partial_1} R_{\alpha}.
\end{equation}
where $\alpha \in \Cl(X_{\Sigma})$. It is known that $K_{\bullet}(\f)$ provides a minimal free resolution of $I(\f)$, if and only if, $I(\f)$ is a complete intersection. Moreover, for any specialization $\f$ and any $\alpha \in \Cl(X_{\Sigma})$, we have that $\partial_1$ equals
\begin{equation}
\label{macaulay2}
    \bigoplus_{i = 0}^nR_{\alpha - \alpha_i} \xrightarrow[]{} R_{\alpha}, \quad (\G_0,\dots,\G_n) \xrightarrow[]{} \sum_{i = 0}^n\G_i\F_i.
\end{equation} 
Using \eqref{monomialbasis}, we can write the matrix of this map in the monomial bases given by the polytopes $\Delta_{[D - D_i]}$ for $\alpha = [D]$. Observe that, if $\alpha$ corresponds to the convex hull of $\mathcal{B}$ as in \eqref{banddelta}, the matrix of \eqref{macaulay2} coincides with the matrix $\mathcal{M}_{\mathcal{A}}$ from \eqref{macaulay}.

In \cite{gkz1994}, they show that for $\alpha \gg 0$ (roughly speaking, for sufficiently large polytopes $\Delta_{D}$), the homology of the complex is concentrated in the last term, i.e. $K_{\bullet}(\f)_{\alpha}$ is exact, if and only if,
$$ H_0(K_{\bullet}(\f))_{\alpha} = (R/I(\f))_{\alpha} = 0.$$

If there exists $\alpha \in \operatorname{Cl}(X_{\Sigma})$ such that $(R/I(\f)){\alpha} = 0$, then $V(I(\f)) \subset X{\Sigma}$ is empty. In other words, the non-exactness of the Koszul complex detects precisely the same condition as the resultant: the solvability of the system over $X_{\Sigma}$. This suggests that the resultant can be computed as the determinant of this complex. We present Cayley's formula as the main tool for computing such determinants.

\begin{definition}
    Let $F_{\bullet}$ be a complex of $\mathbb{K}$-vector spaces of length $k$, $r_i = \dim F_i$ for $i = 1,\dots,k$ and $M_i$ be the matrix of $\partial_i$ in some basis. Consider a tuple of subsets $(K_0,\dots,K_n)$ such that $I_i \subset [r_i]$, $K_0 = [r_0]$ and the submatrix $\mathcal{H}_i$ of $M_i$ indexed by $B_i \setminus K_i$ and $K_{i-1}$ is square and invertible. The determinant of $F_{\bullet}$ is
    $$\det(F_{\bullet}) \coloneqq \prod_{i = 1}^k \det(\mathcal{H}_i)^{(-1)^{i}}.$$
    
\end{definition}

It follows from \cite[Appendix A]{gkz1994} that, in analogy to the case of matrices, the determinant of a complex vanishes if the complex fails to be exact and the formula is indepedent of the chosen subsets $(K_0,\dots,K_n)$.

\begin{theorem}
\label{main}

Let $K_{\bullet}(\f)$ be the Koszul complex, as in \eqref{koszulcomplex} for some specialization $\f$ of the generic homogeneous system defined by $D_0,\dots,D_n$. There exists $\alpha \in \Cl(X_{\Sigma})$ such that for every nef divisor class $\alpha' \in \Nef(X_{\Sigma})$, the following are equivalent:
\begin{itemize}
    \item $V(I(\f)) \subset X_{\Sigma}$ is empty.
    \item The last differential in $K_{\bullet}(\f)_{\alpha + \alpha'}$ has full rank.
    \item $K_{\bullet}(\f)_{\alpha + \alpha'}$ is exact.
\end{itemize}
\end{theorem}

The proof can be found in \cite[Theorem 1.7]{sparsenullstellensatz} and in \cite{gkz1994}. Therefore, checking the degree of the determinant of $K_{\bullet}(\F)_{\alpha}$ for $\alpha$ satisfying Theorem \ref{main} suffices in order to derive that
\begin{equation}
\label{determinantkoszul}
\det(K_{\bullet}(\F)_{\alpha}) = \Res_{\mathcal{A}}.
\end{equation} The proof can be found in \cite[Chapter 3, Theorem 4.2]{gkz1994}, but intuition follows from noting that the degree $\det(K_{\bullet}(\F)_{\bullet})$ on the coefficients of $\F_i$ is an alternate sum of the dimensions of $R_{\alpha}$, which count the number of lattice points on some polytopes while the degrees of $\Res_{\mathcal{A}}$ equal the mixed volumes $\MV_{-i}(\Delta)$. A relation between the two can be found in \cite{ehrhart}.

Observe the proximity of the discussion about resultants and the sparse nullstellensatz: given polynomials $\f_1,\dots,\f_r$ such that $V(I(\f)) \subset X_{\Sigma}$ is empty, find the smallest $\alpha$ such that $(R/I(\f))_{\alpha} = 0$.  While for affine polynomials, the bounds can be exponential in the number of variables \cite{10.1006/aama.1998.0633}, when we want to impose that there are no solutions over $X_{\Sigma}$ this bound (on the degree) is always linear in the number of variables and has appeared in \cite{tuitman,sparsenullstellensatz}.

In general, the goal is to find degrees $\alpha \in \Cl(X)$ such that $\Res_{\mathcal{A}}$ is the determinant of $K_{\bullet}(\F)_{\alpha}$. The following theorems on vanishing of sheaf cohomology over toric varieties are useful to decide when this happens. Proofs of such results are given in \cite[Theorem 9.2.3, Theorem 9.2.7, Theorem 9.3.5]{coxlittleschneck}.
\begin{theorem}
\label{vanishing}
Let $D = \sum_{j = 1}^{n + r} a_jD_j$ be a Weil divisor, as above. Then,
\begin{itemize}
    \item (Demazure) If $D$ is nef, then
    $H^i(X_{\Sigma},\mathcal{O}(D)) = 0$ for all $i > 0$.
    \item (Batyrev-Borisov) If $D$ is nef, then
$H^{n-i}(X_{\Sigma},\mathcal{O}(-D)) = 0$ for all $i > 0$.
    \item (Mustata) If $D$ is a $\mathbb{Q}$-Cartier nef divisor, i.e. there exists $c \in \mathbb{Z}_{>0}$ such that $cD$ is nef and Cartier, then $H^i(X_{\Sigma},\mathcal{O}(\lfloor D \rfloor)) = 0$ for $i > 0$ where $\lfloor D \rfloor = \sum \lfloor a_j \rfloor D_j$.
    \end{itemize}
\end{theorem}

As all of these statements relate to nef divisors, if $D,D'$ are both nef, then the statment also holds for the divisor $D + D'$. This justifies the statement of Theorem \ref{main} because if the equivalences hold for $\alpha$, they must also hold for $\alpha + \alpha'$.

The third item in Theorem \ref{vanishing} has been used recently by D'Andrea and Jerónimo \cite{sparsenullstellensatz} to show that Theorem \ref{main} holds for $\alpha = [D_{\delta}]$ defined as $D_{\delta} = \sum_{i = 0}^n D_i - D_\delta$ for $\delta$ as in \eqref{banddelta} and $D_{\delta}$ defined as $\sum_{j = 1}^{n+ r}c_jD_j$ where $c_j = 1$ if $\langle u_j,\delta \rangle > 0$ and $0$ otherwise. This is exactly the Newton polytope of the lattice points in translated cells $\mathcal{B}$, which was used in the construction of the Canny-Emiris formula. Note that it is possible that the Newton polytope of $\mathcal{B}$ does not correspond to a nef divisor in $X_{\Sigma}$. For this reason, one can consider a common refinement $\Sigma'$ of the normal fans of $\Delta_i$ and the convex hull of $\mathcal{B}$.

When $X_{\Sigma}$ is a smooth toric variety, one can say a bit more: the Hilbert function of $I(\f)$ becomes a polynomial for $\alpha \gg 0$ which can be expressed using local cohomology modules \cite[Lemma 2.8, Proposition 2.14]{boundsmsmith}. Classically, in the $0$-dimensional case, the Hilbert polynomial counts the number of solutions with multiplicity. Therefore, the same vanishing results can be used (also beyond the case of $n + 1$ polynomials) to show that if $V(I(\f))$ is zero-dimensional, then the corank of $\partial_1$ equals the number of points in  $V(I(\f))$, counted with multiplicity \cite{BUSE2024107739, BUSE2022514}. Even beyond the smooth case, similar results can also be obtained \cite{bender2021toric}.

\section{Resultant polytopes}
\label{sectionpolytopes}

In this section, we review a different approach towards sparse resultants. Instead of seeking an exact determinantal representation, we focus on its Newton polytope $N(\Res_{\mathcal{A}})$, also known as \textit{resultant polytope}, or an   orthogonal projection of $N(\ResA)$ along a given direction.
These projections are related to specializations of the resultant and can be used in implicitization, system solving using the $u$-resultant and other applications. Given its Newton polytope, one can compute $\Res_{\mathcal{A}}$ (or the eliminant, in some cases) by interpolating some of its values.

Methods for computing the resultant polytope are strongly combinatorial and have been approached using fiber polytopes \cite{billerasturmfels, sturmfels94}, mixed fiber polytopes \cite{esterov2010newton, estkho,mohrmukhina} and tropical geometry \cite{StuYu08, JensenYu11, sturmfelsyu}.
In this survey, we focus on an algorithm that was presented in \cite{oraclebased, EmFiKoPeJ} which, using an oracle, produces  vertices of the resultant polytope, or its orthogonal projection, in a
given direction, and avoids walking on the polytope whose dimension increases with the number of lattice points in the tuple $\mathcal{A} = (\mathcal{A}_0,\dots,\mathcal{A}_n)$. The approach is output-sensitive as it makes one oracle call per vertex and facet. It extends to other polytopes defined in elimination theory, such as \textit{secondary polytopes} and Newton polytopes of discriminants. Computing projections of $N(\Res_{\mathcal{A}})$ is relevant because often (for instance, in polynomial system solving or implicitization) most of the coefficients of the resultant are specialized, except for a small number, see Example ~\ref{CayleyExample}.

Although tropical geometry exceeds the scope of this survey, many of the methods to recover the resultant polytope, the Newton polytope of the implicit equation of parametric curves or surfaces, or the Newton polytope of the discriminant, extensively use tropical techniques \cite{JensenYu11,StuYu08, tropres, tropicalimplicitization, cueto2012implicitizationsurfacesgeometrictropicalization,Rincon12}.

\begin{figure}[t]
    \centering
    \includegraphics[width=0.75\textwidth]{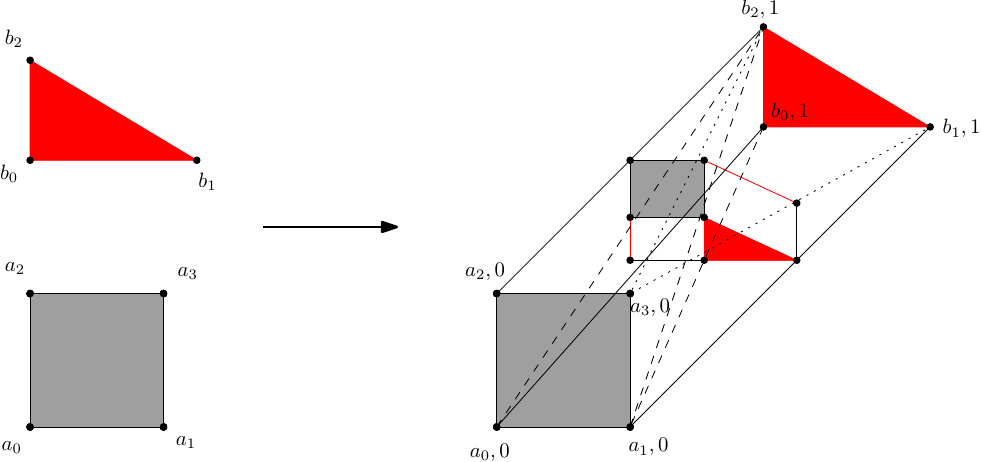}
    \caption{The Cayley trick applied to a triangle and a square.}
    \label{fig:CayleyTrick}
\end{figure}

Given the supports of the input polynomials, we first define an auxiliary set.

\begin{definition}
Given the family $\mathcal{A} = (\mathcal{A}_0,\dots,\mathcal{A}_n)$, $\mathcal{A}_i\subset \mathbb{Z}^n$, we define the 
\textbf{Cayley set} as
\begin{equation}\label{EQ:Cayley}
\Cay(\mathcal{A}):=\bigcup_{i=0}^{n} (\mathcal{A}_{i} \times \{e_{i}\}) \subset \ZZ^{2n},
\end{equation}
where $e_0,\ldots,e_n \in \mathbb{Z}^n$ form an affine basis of $\RR^n$, i.e.,
$e_{0}$ is the zero vector, $e_i=(0,\dots,0,1,0,\dots,0)$, $i=1,\dots,n$.
Clearly, $|\Cay(\mathcal{A})| = \sum_i |\mathcal{A}_i|$, where $|\cdot|$ denotes cardinality. We denote by $M_{\mathcal{A}}$ the matrix whose columns are the supports in the Cayley set.
\end{definition}

The following theorem establishes a correspondence between the triangulations of $\Cay(\mathcal{A})$ and mixed subdivisions of $\Delta=\sum_i\conv(\Ac_i)$, and it is known as the \textit{Cayley trick} \cite{gkz1994}, see also \Cref{fig:CayleyTrick}.

\begin{proposition}\label{P:Cayley_trick}
There exists a bijection between
coherent tight mixed subdivisions of $\Delta$ and
coherent triangulations of $\Cay(\Ac)$. In particular, there exists a bijection between tight mixed subdivisions of $\Delta$ and triangulations of $\Cay(\Ac)$
the mixed subdivisions of  $\Delta$ and the polyhedral subdivisions $\Cay(\Ac)$. 
\end{proposition}

\begin{example}
\label{CayleyExample}
A standard benchmark in geometric modeling is the implicitization of the
bicubic surface defined by three polynomials in two parameters and total degrees $3,3,6$, respectively.
The supports $\mathcal{A}_0,\mathcal{A}_1,\mathcal{A}_2 \subset \mathbb{Z}^2$ of the input polynomials  
have cardinalities $7,6,14$. 
The Cayley set $\Cay(\Ac)\subset\ZZ^4$, constructed as in \eqref{EQ:Cayley}, has
$7+6+14=27$ points. It is depicted in the following matrix, with coordinates as
columns, where the supports from different polynomials and the Cayley
coordinates are distinguished.
%

\setlength{\tabcolsep}{0pt} 
\newcommand{\w}{4mm} 
\[M_{\mathcal{A}}=
\begin{tabular}{|m{.5mm} >{\centering}m{\w} >{\centering}m{\w}
>{\centering}m{\w} >{\centering}m{\w} >{\centering}m{\w} >{\centering}m{\w}
>{\centering}m{\w}| 
>{\centering}m{\w} >{\centering}m{\w} >{\centering}m{\w}
>{\centering}m{\w} >{\centering}m{\w} >{\centering}m{\w}|
>{\centering}m{\w}
>{\centering}m{\w} >{\centering}m{\w} >{\centering}m{\w} >{\centering}m{\w}
>{\centering}m{\w} >{\centering}m{\w} >{\centering}m{\w} >{\centering}m{\w}
>{\centering}m{\w} >{\centering}m{\w} >{\centering}m{\w} >{\centering}m{\w}
>{\centering}m{\w} m{.5mm}| l}

\cline{1-1}
\cline{29-29}

& \color{Maroon} $0$ & \color{Maroon}  $0$ & \color{Maroon}  $1$ &
\color{Maroon} $0$ & \color{Maroon} $2$ & \color{Maroon} $0$ & \color{Maroon}
$3$ &

\color{Blue} $0$ & \color{Blue} $0$ & \color{Blue} $1$ & \color{Blue}  $2$ &
\color{Blue}  $0$ & \color{Blue} $3$ &

 \color{ForestGreen} $0$ & \color{ForestGreen} $0$ & \color{ForestGreen} $1$ &
\color{ForestGreen} $0$ &
\color{ForestGreen} $1$ & \color{ForestGreen} $2$ & \color{ForestGreen}  $1$ &
\color{ForestGreen}  $2$ &
\color{ForestGreen} $1$ & \color{ForestGreen} $2$ & \color{ForestGreen} $3$ &
\color{ForestGreen} $2$ &
\color{ForestGreen} $3$ & \color{ForestGreen} $3$ & \color{ForestGreen} &
        \multirow{2}{*}{ \Big\} \text{support}}\\

& \color{Maroon} $0$ & \color{Maroon}  $1$ & \color{Maroon}  $0$ &
\color{Maroon} $2$ &
\color{Maroon} $0$ & \color{Maroon} $3$ & \color{Maroon} $0$ & 

\color{Blue} $0$ & \color{Blue} $1$ & \color{Blue} $0$ & \color{Blue}  $0$ &
\color{Blue}  $3$ & \color{Blue} $0$ &

 \color{ForestGreen} $0$ & \color{ForestGreen} $1$ &
\color{ForestGreen} $0$ & \color{ForestGreen} $2$ & \color{ForestGreen} $1$ 
& \color{ForestGreen} $0$ & \color{ForestGreen}  $2$ & \color{ForestGreen}  $1$
& \color{ForestGreen} $3$ &
\color{ForestGreen} $2$ & \color{ForestGreen} $1$ & \color{ForestGreen} $3$ &
\color{ForestGreen} $2$ &
\color{ForestGreen} $3$ & \color{ForestGreen} &\\

\cline{2-28}

& \color{BurntOrange} $0$ & \color{BurntOrange}  $0$ & \color{BurntOrange}  $0$
& \color{BurntOrange} $0$ &
\color{BurntOrange} $0$ & \color{BurntOrange} $0$ & \color{BurntOrange} $0$ &
\color{BurntOrange} $1$ &
\color{BurntOrange} $1$ 
& \color{BurntOrange} $1$ & \color{BurntOrange}  $1$ & \color{BurntOrange}  $1$
& \color{BurntOrange} $1$ &
\color{BurntOrange} $0$ & \color{BurntOrange} $0$ & \color{BurntOrange} $0$ &
\color{BurntOrange} $0$ &
\color{BurntOrange} $0$ 
& \color{BurntOrange} $0$ & \color{BurntOrange}  $0$ & \color{BurntOrange}  $0$
& \color{BurntOrange} $0$ &
\color{BurntOrange} $0$ & \color{BurntOrange} $0$ & \color{BurntOrange} $0$ &
\color{BurntOrange} $0$ &
\color{BurntOrange} $0$ & \color{BurntOrange} &
        \multirow{2}{*}{ \Big\} \text{Cayley}}\\

& \color{BurntOrange} $0$ & \color{BurntOrange}  $0$ & \color{BurntOrange}  $0$
& \color{BurntOrange} $0$ &
\color{BurntOrange} $0$ & \color{BurntOrange} $0$ & \color{BurntOrange} $0$ &
\color{BurntOrange} $0$ &
\color{BurntOrange} $0$ 
& \color{BurntOrange} $0$ & \color{BurntOrange}  $0$ & \color{BurntOrange}  $0$
& \color{BurntOrange} $0$ &
\color{BurntOrange} $1$ & \color{BurntOrange} $1$ & \color{BurntOrange} $1$ &
\color{BurntOrange} $1$ &
\color{BurntOrange} $1$ 
& \color{BurntOrange} $1$ & \color{BurntOrange}  $1$ & \color{BurntOrange}  $1$
& \color{BurntOrange} $1$ &
\color{BurntOrange} $1$ & \color{BurntOrange} $1$ & \color{BurntOrange} $1$ &
\color{BurntOrange} $1$ &
\color{BurntOrange} $1$ & \color{BurntOrange} &\\

\cline{1-1}
\cline{29-29}

\end{tabular}
\] 
By \cite[Theorem 2E.1]{gkz1994} 
it follows that $N(\ResA)$ 
has dimension $|\Ac|-4-1=22$ and lies in $\RR^{27}$.
Implicitization requires eliminating the two parameters to obtain a
constraint equation over the symbolic coefficients of the polynomials.
Most of the symbolic coefficients are specialized except for $3$
variables, hence the sought for implicit equation of the surface 
is trivariate and the projection of $N(\ResA)$ lies in $\RR^3$.
The number of triangulations of $\Cay(\Ac)$ is 6 orders of magnitude greater
than the number of vertices of $N(\ResA)$.
The 3-dimensional projection of $N(\ResA)$, which equals the
Newton polytope of the implicit equation, is defined by six vertices.
\end{example}

Let $P\subset\RR^d$ be a pointset whose convex hull is of dimension $d$. 
For any triangulation $\mathcal{T}$ of $P$, define vector $\phi_{\mathcal{T}}\in\RR^{|P|}$
with coordinates
\begin{equation}\label{Evolume_embed}
\phi_{\mathcal{T}}(a)= \sum_{\sigma\in \mathcal{T} : a\in \sigma} \mbox{vol}(\sigma), \qquad a\in P,
\end{equation}
summing over all simplices $\sigma$ of $\mathcal{T}$ having $a$ as a vertex; 
the secondary polytope $\Secpol(P)$
is the convex hull of $\phi_{\mathcal{T}}$ for all triangulations $\mathcal{T}$. The following proposition establishes the dimension and key properties of the secondary polytope.
\begin{proposition}
\label{secondarypolytope}
The vertices of $\Secpol(P)$ correspond to the coherent triangulations of $P$,
while its face lattice corresponds to the poset 
of coherent polyhedral subdivisions of $P$, ordered by refinement. 
A lifting vector produces a coherent triangulation $T$
(resp.\ a coherent polyhedral subdivision of $P$)
if and only if it lies in the normal cone of vertex $\phi_T$
(resp.\ of the corresponding face) of $\Secpol(P)$.  
The dimension of $\Secpol(P)$ is $|P|-d-1$. 
\end{proposition}

At this point, the Newton polytope of the resultant comes into play since, as we have already seen in \Cref{section:cannyEmiris}, we can use a generic mixed subdivision $S(\rho)$ to produce vertices of the resultant polytope. The following result, first proved in \cite[Theorem 2.1]{sturmfels94}, can also be seen as special case of \Cref{initialFormTheorem} which was proven in \cite{dandrea2020cannyemiris}.

\begin{proposition} \label{PSturmf_extreme}
For a sufficiently generic lifting function $\rho$, given by values $\omega \in \prod_{i = 0}^n\mathbb{R}^{\mathcal{A}_i}$, and providing a mixed subdivision $S(\rho)$,
we have
\begin{equation}\label{Eq:Sturmf_extreme}
\init_{\omega}(\Res_{\mathcal{A}}(u_{i,a}t^{\omega_{i,b}})) = \pm \prod_{i=0}^{n} \prod_{\mathcal{D}} u_{i,b}^{\mathrm{vol}(\sigma)},
\end{equation}
where $\mathcal{D}$ ranges over all $n$-dimensional $i$-mixed cells of $S(\rho)$ and $b \in \mathcal{A}_i$ indexes a unique $0$-dimensional component in $\mathcal{D}$.
\end{proposition}

Once this is established, one can show \cite[Theorem 2E.1]{gkz1994}, \cite[Corollary 5.1]{sturmfels94} that $N(\Res_{\mathcal{A}})$ is a Minkowski summand of 
$\Secpol(\Cay(\Ac))$, and that both $\Secpol(\Cay(\Ac))$ and $N(\ResA)$ have 
dimension $|\Ac| -2n-1$.
In  \cite{michielscools} they describe all Minkowksi summands of the secondary polytope.
A linear space $H_{\mathcal{A}}$ of codimension $2n + 1$ in $\RR^{\mathcal{A}}$ in which $N(\Res_{\mathcal{A}})$ is contained can be described by the homogeneities of the sparse resultant \cite[Chapter 7, Proposition 1.11]{gkz1994}. Note that if the family $\mathcal{A}$ is essential, then the projection along $2n + 1$ affinely independent points of $\Cay(\mathcal{A})$ produces an isomorphic polytope.

To avoid the high complexity of directly computing $\Secpol(\Cay(\Ac))$, which has many more vertices than $N(\Res_{\mathcal{A}})$ (Example ~\ref{CayleyExample}), one can instead define an equivalence relation on its vertices. In \cite{MicVer99}, they define an equivalence class on the vertices of the secondary polytope for which the corresponding mixed subdivisons have the same mixed cells. The classes map in a many-to-one fashion to resultant vertices.
The algorithm in \cite{EmFiKoPeJ} exploits a stronger equivalence relation to design an oracle that, given a direction, returns the corresponding resultant vertex. 

\begin{figure}[t] \centering
 \includegraphics[width=1\textwidth]{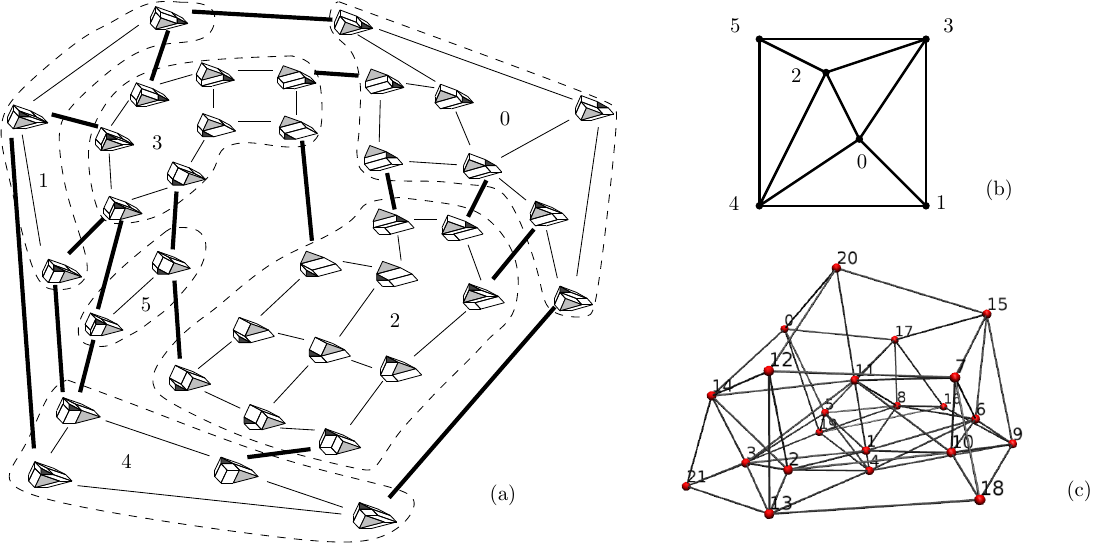}
\caption{(a) 
The secondary polytope $\Secpol(\Cay(\Ac))$ of two triangles (dark, light grey) and one
segment $\Ac_0=\{(0,0),(1,2),(4,1)\},\,\Ac_1=\{(0,1),(1,0)\},\,\Ac_2=\{(0,0),(0,1),(2,0)\}$,
where $\Cay(\Ac)$ is defined as in \eqref{EQ:Cayley};
its vertices correspond to mixed subdivisions of the
Minkowski sum $\Ac_0+\Ac_1+\Ac_2$ and edges to flips between them
(b)
$N(\ResA)$, whose vertices correspond to the dashed classes of $\Secpol(\Cay(\Ac))$.
Bold edges of $\Secpol(\Cay(\Ac))$, called cubical flips, map to edges of $N(\Res_{\mathcal{A}})$
(c)
$4$-dimensional $N(\Res_{\mathcal{A}})$ of 3 generic trinomials.
}
\label{fig:sec_res}
\end{figure}

\begin{figure}[t]
\centering
\includegraphics[scale=0.9]{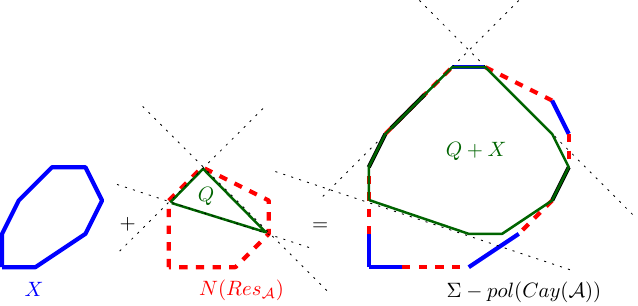}
\caption{The key idea of the algorithm: each illegal hyperplane of $Q$ with normal
$w$, separates
the already computed vertices of $\varPi$ (here equal to $N(\Res_{\mathcal{A}})$) from new ones,
extremal with respect to $w$. In the figure, $X$ is a polytope
such that $X+N(\Res_{\mathcal{A}})=\Secpol(\Cay(\Ac))$.}
\label{fig:oneperhplane}
\end{figure}

In some cases, it is possible to theoretically characterize $N(\Res_\mathcal{A})$:
\begin{itemize}
    \item In~\cite{GELFAND1990237}, Gelfand, Kapranov and Zelevinsky show that $N(\Res_{\mathcal{A}})$,
for two univariate polynomials with $k_0+1,k_1+1$
monomials, has $\binom{k_0+k_1}{k_0}$ vertices and, when both $k_i\ge 2$,
it has $k_0k_1+3$ facets.

\item In \cite{sturmfels94}, Sturmfels showed that $N(\Res_{\mathcal{A}})$ is 1-dimensional
if and only if $|A_i|=2$, for all $i$, the only planar $N(\Res_{\mathcal{A}})$
is the triangle, whereas the only 3-dimensional ones are the tetrahedron,
the square-based pyramid, and the resultant polytope of two univariate trinomials. In \Cref{fig:sec_res}(b), a polytope isomorphic to the resultant polytope of three bivariate polynomials is shown.

\item Following \cite[Theorem~6.2]{sturmfels94}, the 4-dimensional polytopes
include the 4-simplex, some $N(\Res_{\mathcal{A}})$ obtained by pairs of univariate polynomials,
and those of~3 trinomials, which have been investigated in~\cite{DEF12} using the algorithm presented here. 
Among the latter, the polytope with the largest number of vertices that was computed using the algorithm in \cite{oraclebased} appears in \Cref{fig:sec_res}(c); its f-vector is $(22,66,66,22)$. 
\end{itemize}

The main contribution of \cite{EmFiKoPeJ} was an oracle-based algorithm for computing $N(\Res_{\mathcal{A}})$ or some projection of it, which can be viewed as the Newton polytope of generic specializations of the sparse resultant.
The algorithm utilizes the \textit{beneath and beyond} method to compute 
both vertex and half-space representations. Its incremental nature implies that we also obtain a triangulation of the
polytope, which may be useful to enumerate its lattice points. 

The algorithm also extends to computing $\Secpol(\Cay(\Ac))$, as well as 
the Newton polytope of the discriminant and, more generally, any polytope that can
be efficiently described by a vertex oracle or its orthogonal projection. For example, it suffices to replace the oracle by the one in~\cite{Rincon12} to obtain a method for computing the discriminant polytope. There exists a publicly available implementation of the algorithm based on CGAL~\cite{CGAL}.

To compute an orthogonal projection of $N(\Res_{\mathcal{A}})$, denoted as
$\varPi \subset \RR^m$, one considers the map:
$$
\pi:\RR^{|\Ac|}\rightarrow \RR^m : N(\ResA)\rightarrow \varPi,\;\, m \le |\Ac| .
$$
By reindexing, this is the subspace of the first $m$ coordinates,
so $\pi(w)=(w_1,\dots,w_m)$.
It is possible that none of the coefficients $u_{ij}$ is specialized,
hence $m = |\Ac|$, $\pi$ is trivial, and $\varPi=N(\ResA)$.
Assuming that the specialized coefficients take sufficiently
generic values, $\varPi$ is the Newton polytope of the corresponding
specialization of $\ResA$.
The following result, due to Jensen and Yu \cite[Lemma~3.20]{JensenYu11}, simplifies the computation of the projection.

\begin{lemma}\label{Linsidepoints}
If $b_i \in \Ac_i$ corresponds to a specialized coefficient of $F_i$
and lies in the convex hull of the other points in $\Ac_i$ corresponding
to specialized coefficients, then removing $a_{ij}$ from $\Ac_i$ does not
change the Newton polytope of the specialized resultant.
\end{lemma} 

%

\begin{definition}
A hyperplane $H$ is called \emph{legal}
if it is a supporting hyperplane to a facet of $\varPi$,
otherwise it is called \emph{illegal}.
\end{definition}

In order to compute $\Pi$, the algorithm uses Propositions  \ref{P:Cayley_trick} and  \ref{PSturmf_extreme} to compute vertices of the resultant polytope.
Initially, it computes a polytope $Q \subset \varPi$ such that $\dim(Q)=\dim(P)$, by sampling a list of general directions in the linear space $\pi(H_{\mathcal{A}})$. Then, using the facet representation of $Q$ one can list the supporting hyperplanes and use the following .

\begin{lemma}\cite[Lemma 5]{EmFiKoPeJ}
Let $v \in N(\mathcal{R})$ be a vertex produced by a normal vector $w$ to a supporting hyperplane $H$ of $Q$.
Then $v\not\in H$ if and only if $H$ is not a supporting hyperplane of $\varPi$.
\end{lemma}

If $v\notin H$, it is a new vertex thus yielding a tighter \textit{inner
approximation} of $\varPi$ by inserting it to $Q$, i.e.\ $Q\subset
\conv(Q\cup v) \subseteq \varPi$. 
This happens when the preimage  $\pi^{-1}(\mathcal{F})\subset N(\Res_{\mathcal{A}})$ of the facet $\mathcal{F}$ of
$Q$ defined by $H$, is not a Minkowski summand of a face of $\Sigma(\A)$ having normal $w$ (\cref{fig:oneperhplane}).
Otherwise, there are two cases: either $v\in H$ and $v\in Q$, 
thus the algorithm simply decides that hyperplane $H$ is legal, or
$v\in H$ and $v \notin Q$, in which case the algorithm again
decides that $H$ is legal but also inserts $v$ to $Q$.

As a consequence, one can enumerate all supporting hyperplanes of $Q$ and distinguish those that are legal or ilegal, by computing the vertex corresponding to the outer normal of that hyperplane. The algorithm terminates with both the vertex and facet representations of $\varPi$ and makes only one call for each vertex and each facet of the resultant polytope.

\printbibliography

\newcommand{\Adresses}{%
    \bigskip
    \footnotesize
    C. Checa, \textsc{Department of Mathematics, University of Copenhagen, Universitetsparken 5, 2100, Denmark}, \texttt{ccn@math.ku.dk}
    \medskip
    
    I.~Z. Emiris, \textsc{Department of Informatics \& Telecommunications, University of Athens, and Heron Center of Excellence in Robotics, Athena Research Center, 6 Artemidos \& Epidavrou 15125, Greece}, \texttt{emiris@athenarc.gr}
    \medskip
    
    C. Konaxis, \textsc{Athena Research Center, 6 Artemidos \& Epidavrou 15125, Greece},  \texttt{ckonaxis@athenarc.gr}
}

\Adresses

\end{document}